\documentclass[letterpaper, 10 pt, conference]{ieeeconf}  

\IEEEoverridecommandlockouts                              
\usepackage{graphics} 
\usepackage{epsfig} 
\usepackage{mathptmx} 
\usepackage{times} 
\usepackage{amsmath} 
\usepackage{amssymb}  
\usepackage{xcolor}
\usepackage{booktabs}
\usepackage{caption}
\usepackage{tabulary}
\usepackage{tabularx}
\usepackage{cprotect}

\usepackage{algorithm}
\usepackage[noend]{algpseudocode}

\makeatletter
\def\BState{\State\hskip-\ALG@thistlm}
\makeatother

\allowdisplaybreaks

\newcommand*\xbar[1]{%
  \hbox{%
    \vbox{%
      \hrule height 0.5pt 
      \kern0.5ex
      \hbox{%
        \kern-0.1em
        \ensuremath{#1}%
        \kern-0.1em
      }%
    }%
  }%
} 

\newtheorem{rem}{Remark}

\title{\LARGE \bf
Estimating effective connectivity in linear brain network models}

\author{G. Prando, M. Zorzi, A. Bertoldo and A. Chiuso $^\dagger$
\thanks{This work has been partially supported by the  project ``Statistical learning methods for estimating the effective connectivity of human brain networks'' (BIRD162411/11).}
\thanks{$^\dagger$ Dept. of Information  Engineering, University of Padova (e-mail: \{\tt \footnotesize prandogi,zorzimat,bertoldo,chiuso\}@dei.unipd.it)}%
}

\begin{document}

\maketitle
\thispagestyle{empty}
\pagestyle{empty}

\begin{abstract}
Contemporary neuroscience has embraced network science to study the complex and self-organized structure of the human brain;  one of the main outstanding issues is that of inferring from measure data, chiefly functional Magnetic Resonance Imaging (fMRI),  the so-called effective connectivity in brain networks, that is the existing interactions among neuronal populations. This inverse problem is complicated  by the fact that the BOLD (Blood Oxygenation Level Dependent) signal measured by  fMRI  represent a dynamic and nonlinear transformation (the hemodynamic response) of neuronal activity.
In this paper, we consider resting state (rs) fMRI data; building upon a linear population model of the BOLD signal and a stochastic linear DCM model, the model parameters are estimated through an EM-type iterative procedure, which alternately estimates the neuronal activity by means of the Rauch-Tung-Striebel (RTS) smoother, updates the connections among neuronal states and refines the parameters of the hemodynamic model; sparsity in the interconnection structure is favoured using an iteratively reweighting scheme. 

Experimental results using rs-fMRI data are shown demonstrating the effectiveness of our approach and comparison with state of the art routines (SPM12 toolbox) is provided.


\end{abstract}

\section{Introduction}
Recent research in the field of computational neuroscience has shown a growing interest in the estimation of the effective connectivity among brain regions, that is, the existing 
interactions between neuronal populations located in different brain areas \cite{friston2011functional}. What makes the problem particularly challenging is the lack of a direct measurement of the neuronal activity: standard neuroimaging techniques (such as EEG, MEG and fMRI) only provide an indirect measurement of such activity. In particular we shall be concerned with  fMRI data, which provide noisy measurements of  the so-called BOLD (Blood Oxygenation Level Dependent) signal.
The dependency between this quantity and neuronal activity is described by the so-called hemodynamic response. In brief, an increment of the synaptic activity in a certain brain region causes an increased blood flow to that area, which in turn augments the ratio of oxygenated haemoglobin w.r.t. the deoxygenated one in that region. By exploiting the different magnetic properties of these two types of hemoglobin, fMRI is able to detect local changes in their relative concentration, thus highlighting the eventual presence of synaptic activations.\\
Existing methods tackle the problem of effective connectivity estimation by postulating a parametric generative model of the observed fMRI data, which not only accounts for the coupling among neuronal populations, but also for the aforementioned hemodynamic response. Nonlinear dynamical models are used and generally referred to as Dynamic Causal Models (DCMs) \cite{friston2003dynamic}: both deterministic \cite{friston2003dynamic,marreiros2008dynamic,stephan2008nonlinear} and stochastic \cite{li2011generalised,friston2014dcm} formulations exist, with the latter being more suited to describe resting state fMRI data.\\
Having postulated a DCM,  effective connectivity inference is reduced to model parameters estimation, which is typically accomplished using Variational Bayes techniques \cite{friston2007variational,friston2008variational,daunizeau2009variational}. However, existing approaches are computationally intensive and are effective only on small brain networks. A large portion of the required computational effort is spent for the determination of the structure of the network connecting the different brain regions. Indeed, this 
choice involves the comparison of several candidate structures, making the problem combinatorial, \cite{friston2011post}. In addition, the inversion of stochastic DCMs appears more involved since it also requires the estimation of the hidden random neuronal activity \cite{daunizeau2009variational}. \\
DCM is a Bayesian approach originally proposed to investigate from fMRI signal the influence on the neuronal system exerts by external stimuli, i.e. the so-called task-evoked effective connectivity. DCM can be also be used to study the directed connections among distributed brain areas when only endogenous brain fluctuations are present, i.e. the so-called resting-state effective connectivity. Our study will focus on this latter problem. We develop a linear DCM, and adapt    filtering techniques to estimate the hidden neuronal activity as well as the DCM parameters. Differently from existing methodologies, which require the comparison of a set of candidate patterns of effective connectivity, our  approach estimates directly the connectivity pattern, by resorting to an iteratively $\ell_1$ reweighed algorithm which induces sparsity in the estimated  connectivity matrix, thus automatically detecting the structure of the brain network. \\
The paper is organized as follows. We briefly review classical DCMs and parameters inference techniques in Sec.~\ref{sec:problem_form}, while in
Sec.~\ref{sec:linear_dcm} we propose an approximated linear stochastic DCM, which is suited for resting state fMRI data. In Sec.~\ref{sec:estimation_measured_x} we deals with the estimation of the DCM parameters (unrealistically) assuming to measure the neuronal activity. We deal with the realistic situation in which only the BOLD signal is available in Sec.~\ref{sec:estimation_unmeasured_x}, where we illustrate the joint estimation of the DCM parameters and of the hidden neuronal activity. Using simulated fMRI data, we evaluate the effectiveness of the proposed method in Sec.~\ref{sec:experiments}, while we provide some concluding remarks in Sec.~\ref{sec:conlusion}.

\section{Problem Formulation}\label{sec:problem_form}
\subsection{Generative model}\label{sec:classical_dcm}

A DCM is a nonlinear multi input-multi output (MIMO) dynamic system, with one output for each monitored brain region. A DCM basically
consists of two components: the first  describes the dynamic coupling among neuronal populations, while the second  maps the neuronal activity to the BOLD signal observed in each brain area, through the hemodynamic response. Accordingly, the system states are grouped in two classes: one contains the $n$ neuronal states (one for each monitored brain region), while the other refers to the biophysical quantities involved in the hemodynamic model.\\
In the seminal paper \cite{friston2003dynamic}, the coupling among the neuronal states $x\in\mathbb{R}^n$ is modelled through a deterministic bilinear interaction:
\begin{equation}\label{equ:bilinear_neuronal_model}
\dot{x}(t) = \left(A + \sum_{j=1}^m u_j(t) B_j\right)x(t) + Cu(t), \qquad t\in\mathbb{R}
\end{equation}
where $u(t)=[u_1(t)\ \cdots\ u_m(t)]^T$, with $u_j(t)$ denoting the $j$-th designed external stimulus.
Effective connectivity in the absence of external solicitations is encoded by the matrix $A\in\mathbb{R}^{n\times n}$, while the bilinear term $B_j\in\mathbb{R}^{n\times n}$ accounts for the change in the neuronal coupling due to the $j$-th input; finally, $C\in\mathbb{R}^{n\times m}$ modules the direct influence of external inputs on the neuronal activity. The system matrices $\theta_c=\{A,B_1, B_2, ..., B_m,C \}$ are unknown parameters which need to be estimated in order to retrieve the effective connectivity among the considered brain areas.\\
Extensions of the bilinear model \eqref{equ:bilinear_neuronal_model} can be found in the literature, see e.g. \cite{marreiros2008dynamic} and \cite{stephan2008nonlinear}.\\
In this paper we will consider the resting state activity, which can be modelled as \cite{friston2014dcm}
\begin{equation}\label{equ:stochastic_bilinear_neuronal_model}
\dot{x}(t) = Ax(t) + Cu(t) + w(t), \quad t\in\mathbb{R}
\end{equation}
where $w(t)$ is a stochastic process which accounts for the unmeasurable random fluctuations that drive the resting state activity. Even if in this situation the exogenous inputs $u(t)$ are typically set to zero, model \eqref{equ:stochastic_bilinear_neuronal_model} can still accommodate the presence of external (non-modulatory) signals. In this paper we assume $u(t)=0$.


\noindent The hemodynamic response is generally charachterized through an extension of the so-called Balloon-Windkessel model, \cite{buxton1998dynamics,friston2000nonlinear},  a SISO dynamic system with the neural activity   $x_i(t)$, i.e. the component of $x(t)$ in position $i$,  as input and the corresponding observed BOLD signal $y_i(t)$ as output. The model is  described by $4-th$ order nonlinear state space model  depending upon the  biophysical parameters $\theta_h=\{\kappa,\gamma,\tau,\alpha,\rho\}$ more specifically:
\begin{align}
\dot{s}_i(t) &= x_i(t) - \kappa_i s_i -\gamma_i (f_i(t)-1), \qquad \qquad t\in\mathbb{R}\label{equ:states_balloon_1}\\
\dot{f}_i(t) &= s_i(t)\\
\tau_i \dot{v}_i(t) &= f_i(t) - v_i^{1/\alpha}(t)\\
\tau_i \dot{q}_i(t) &= \frac{f_i}{\rho_i} \left(1- (1-\rho_i)^{1/f_i(t)}\right) - v_i^{1/\alpha}(t)\frac{q_i(t)}{v_i(t)} \label{equ:states_balloon}\\
y_i(t) &= V_0 \left(k_1 (1-q_i(t)) + k_2 \left(1- \frac{q_i(t)}{v_i(t)}\right) + k_3 (1-v_i(t))\right)\nonumber\\
& \hspace{0.5cm}+e_i(t)\label{equ:output_balloon}
\end{align}
The so called hemodynamic states $\{s_i,f_i,v_i,q_i\}$ are biophysical quantities which are affected by the neuronal activity: $s_i$ denotes the vasodilatatory signal $s_i$, $f_i$ is the blood inflow, $v_i$ and $q_i$ are respectively the blood volume and the deoxyhemoglobin content, $e_i(t)$ is the observation error, which may arise from thermal and physiological causes during data acquisition. The output equation \eqref{equ:output_balloon} depends on the resting blood volume fraction $V_0$ (typically set to 0.4) and on the constants $k_1$, $k_2$ and $k_3$. 
Despite these have found different characterizations in the literature, a general form is derived in \cite{stephan2007comparing}:
\begin{equation}
k_1 = 4.3\vartheta_0 \rho_i T_E, \qquad k_2=\varepsilon r_0 \rho_i T_E, \qquad k_3 =1-\varepsilon .
\end{equation}
Here, $T_E$ is the echo time, while $\vartheta_0=40.3s^{-1}$ denotes the frequency offset at the outher surface of the magnetized vessel for fully deoxygenated blood at 1.5 T; $\varepsilon\approx 0.4$ is the ratio of intra- and extravascualr signal, while $r_0=25s^{-1}$ is the slope of the relation between the intravascular relaxation rate and oxygen saturation.

\noindent Collecting the neuronal and the hemodynamic states into the vector $z\in\mathbb{R}^{5n}$ and the two parameters sets in $\theta$, i.e. $\theta=\{\theta_c,\theta_h\}$, a stochastic DCM can be compactly described as the following nonlinear dynamic system:
\begin{equation}\label{equ:full_stochastic_dcm}
\left\{\begin{array}{rcl}
\dot{z}(t) & = & f(z(t),u(t), w(t),\theta) \\
y(t) &=& g(z(t)) + e(t), \qquad \qquad t\in\mathbb{R}
\end{array}\right. 
\end{equation}
where $y(t)=[y_1(t)\ \cdots\ y_n(t)]^T$ and similarly for $e(t)$.
The practical implementation of (\ref{equ:full_stochastic_dcm}) is as follows \footnote{See \texttt{http://www.fil.ion.ucl.ac.uk/spm/}.}. The simulated BOLD data, equivalently fMRI data, are given by discretizing the Balloon model \eqref{equ:states_balloon_1}-\eqref{equ:output_balloon} as well as the trajectory of the neural activity in (\ref{equ:stochastic_bilinear_neuronal_model}) using the sampling time $T_R=2$s, which is the time often used for the data acquisition.  With some abuse of notation,  $x(k):=x(kT_R)$ and $y(k):=y(kT_R)$, with $k\in\mathbb Z$, be, respectively, the discretized trajectory of $x$ and $y$. By \eqref{equ:stochastic_bilinear_neuronal_model}, with $u(t)=0$, it is not difficult to see that 
\begin{equation}\label{equ:rs_neuronal_model}
x(k+1) = e^{A T_R }x(k) + w(k)
\end{equation}
where $w(k)\sim\mathcal{N} (0, Q)$ and
\begin{equation}\label{equ:Q}
Q=\sigma^2 \int_0^{T_R}e^{A \tau}e^{A^T \tau}\mathrm{d}\tau.
\end{equation}

\subsection{Effective connectivity estimation} \label{sec:existing_methods}
When considering the resting state condition, the effective connectivity estimation problem consists in estimating $\theta_c$, i.e. matrix $A$ (and possibly $C$), using the measured BOLD signal $y(k)$, $k=1\ldots N$. It is worth noting that the estimation of $\theta_c$ requires to estimate $\theta_h$ as well, i.e. the parameters of the Balloon model. This appears as a challenging task, because of the nonlinearity of system \eqref{equ:full_stochastic_dcm}, but also of
the unavailability of the random states $z$, which have to be somehow inferred. Previous works have tackled this problem following two different routes.
\\A first approach, referred to as \textit{spectral DCM}, \cite{friston2014dcm}, estimates the second-order statistics (i.e. cross-spectra) of the hidden states $z$, which are assumed to be stationary. In this way, the inference of the hidden time-series $z$ is avoided.
Furthermore, by parametrizing the spectral density of both process and measurement noise (typically as a power law distribution), the DCM is made deterministic and its parameters are estimated through standard Bayesian approaches \cite{friston2003dynamic,friston2008variational}, postulating Gaussian priors on the parameters $\theta$. Typically, a stability-inducing prior is formulated for the matrix $A$, while empirical priors are assigned to the hemodynamic parameters $\theta_h$. However, the nonlinear nature of system \eqref{equ:full_stochastic_dcm} complicates the analytical computation of the posterior probability density function, as well as of the marginal likelihood. To overcome this bottleneck, a variational Bayes scheme is used, under the mean-field and the Laplace assumptions: namely, the approximated marginal posterior distributions are considered independent and Gaussian (see \cite{friston2007variational,daunizeau2009variational} for details).\\
A second family of approaches still parametrizes the distribution of the process and observation noises in \eqref{equ:full_stochastic_dcm}, but operates in the time-domain, thus estimating not only the model parameters $\theta$ but also the hidden neuronal states $z$. The most common methods are Dynamic Expectation Maximization (DEM) \cite{friston2008variational} and Generalized Filtering (GF) \cite{friston2010generalised,li2011generalised}: even if both adopt the Variational Bayes procedure, DEM uses the mean-field and the Laplace approximations, while GF only exploits the latter.\\
Despite these techniques represent the state-of-the-art for the estimation of brain effective connectivity, computational feasibility limits their application to small networks (typically, no more than 10 regions). In particular, this drawback severely affects DEM and GF \cite{razi2015construct}. The reformulation of the nonlinear stochastic DCM \eqref{equ:full_stochastic_dcm} as a linear stochastic dynamic model may represent a way to address this issue. In this paper we pursue this idea, by deriving an approximated linear DCM according to the procedure detailed in Sec.~\ref{sec:linear_dcm}. This reformulation will allow us to jointly estimate the hidden neuronal states $z$ and the DCM parameters $\theta$.
It is worth noting that a recent contribution in this direction focuses on deterministic DCMs and transforms the neuronal state equation \eqref{equ:stochastic_bilinear_neuronal_model} (with $w(t)=0$) into the frequency domain; then, using a fixed hemodynamic response, the inversion of a deterministic DCM is formulated as a Bayesian linear regression problem, \cite{frassle2017regression}.


\section{Statistical  Linearization of DCMs}\label{sec:linear_dcm}
We model the hemodynamic response as a Finite Impulse Response (FIR) model which takes as input only a neuronal state $x_i(k)$ and as output the BOLD signal $y_i(k)$,
\begin{equation}
y_i(k) = \sum_{l=0}^{s-1} h_l x_i(k-l)
\end{equation}
where the length $s$ of the impulse response $h:=[h_1\ \cdots \ h_s]^T$ is chosen large enough to retain all the relevant dynamics components. A prior on the impulse response $h$ is derived, using statistical linearisation techniques, from the Baloon model \eqref{equ:states_balloon_1}-\eqref{equ:output_balloon} where the parameters $\theta_h=\{\kappa,\gamma,\tau,\alpha,\rho\}$ are described by prior distributions typically used in the literature \cite{friston2003dynamic}. This will lead to linear model of the form $h=H\alpha$ (see eq. \eqref{eq:H} below) where  $\alpha$  will have to be estimated under a Bayesian perspective using the available fMRI data. 
To to so, we first  compute a population of typical responses generated by the Balloon model \eqref{equ:states_balloon_1}-\eqref{equ:output_balloon}, then we define $h$ as the linear combination of their empirical mean and of the first $p$ principal components of their sample covariance matrix. 
\noindent Specifically:
\begin{enumerate}
\item Sample $\theta_h^{(j)}$, $j=1,...,N_s$ from the empirical Gaussian distributions reported in Table 1 of the seminal paper \cite{friston2003dynamic}.
\item For each $\theta_h^{(j)}$ compute the output of the Balloon model \eqref{equ:states_balloon_1}-\eqref{equ:output_balloon}, where $x_i(k)=\delta(k)$ and $\delta(k)$ is the Kronecker delta function.  Let  $h^{(j)}$ be the corresponding output truncated at length $s$, which represents a statistic of $h$.
\item Compute the empirical mean $\bar{h}= \frac{1}{N_s}\sum_{j=1}^{N_s} h^{(j)}$.
\item Compute the empirical covariance matrix $\Sigma_h\in\mathbb{R}^{s\times s}$:
\begin{equation}\label{equ:sigma_h}
\Sigma_h  = \frac{1}{N_s}\sum_{j=1}^{N_s} (h^{(j)}-\bar{h})(h^{(j)}-\bar{h})^T.
\end{equation}
\item Compute the eigenvalues decomposition of $\Sigma_h$, 
$\Sigma_h=USU^T$
, where $S:=\mathrm{diag}(s_1,...,s_s)$ and $U:=[u_1\ u_2\ \cdots \ u_s]$.
\item Define $h$ as
\begin{equation} \label{eq:H}
h := H\alpha = \begin{bmatrix}\bar{h} & u_1 & u_2 &\cdots & u_p\end{bmatrix}\alpha, \quad \alpha\in\mathbb{R}^{p+1}.
\end{equation}
\end{enumerate}

\begin{figure}
\centering
\includegraphics[width=\columnwidth]{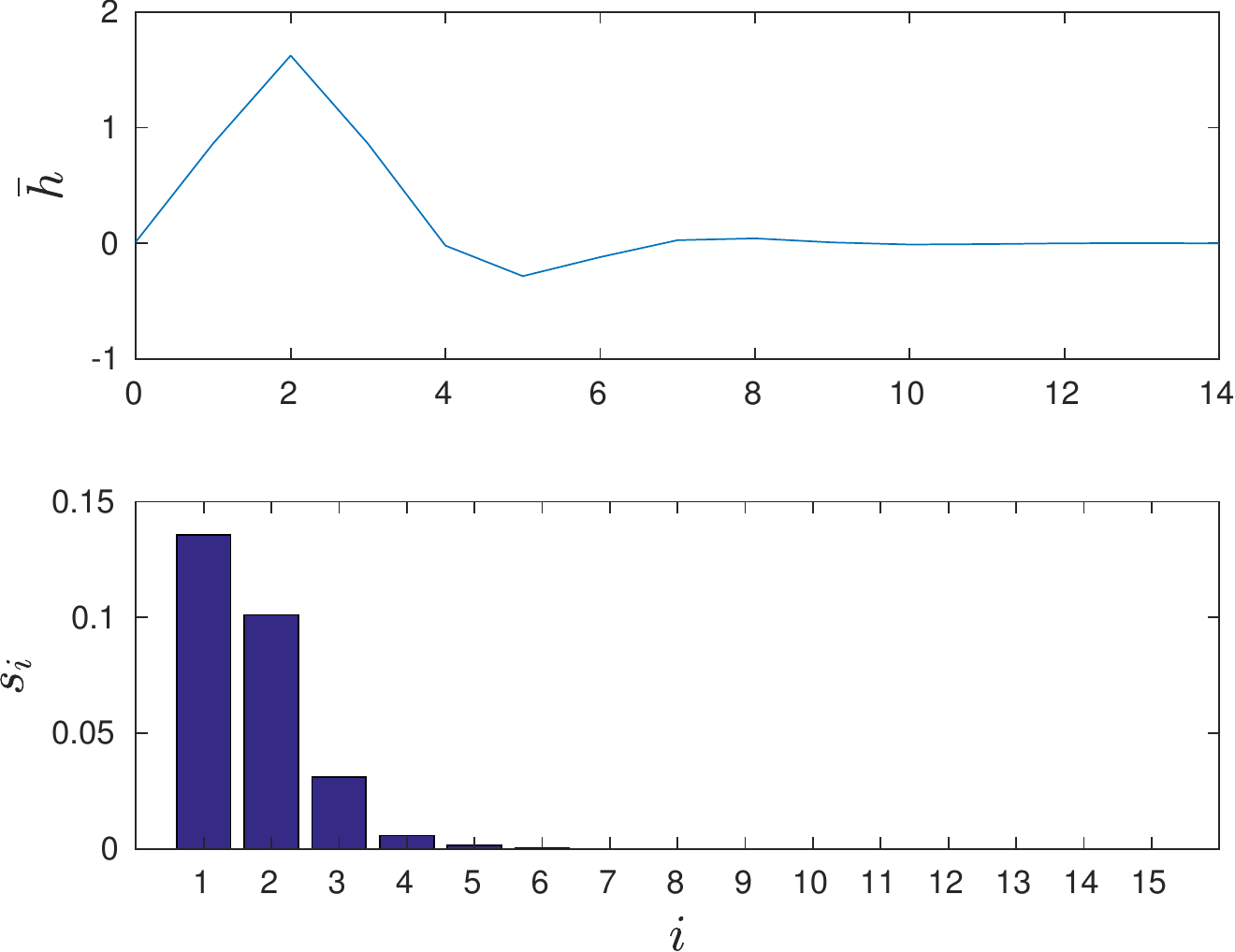}
\caption{\textit{Top:} Empirical mean $\bar{h}$ of the sample hemodynamic responses. \textit{Bottom:} Singular values $s_i$, $i=1,...,s$, of the sample covariance matrix $\Sigma_h$.}\label{fig:hemodynamic}
\end{figure}

\noindent Figure \ref{fig:hemodynamic} illustrates the empirical mean $\bar{h}$ of the sampled population of hemodynamic responses and the singular values of their sample covariance matrix $\Sigma_h$, showing that there is limited variability in the generated responses.

\noindent By introducing the extended state $\mathbf{x}(k)\in\mathbb{R}^{ns}$ 
\begin{align}
\mathbf{x}(k)&:= \begin{bmatrix}
x^T(k) & x^T(k-1) & \cdots & x^T(k-s+1)
\end{bmatrix}^T,
\end{align} 
the matrices
\begin{equation}\label{equ:a_c}
\mathbf{A} :=  \begin{bmatrix}
e^{AT_R} &  \mathbf{0}\\ I_{n(s-1)} & \mathbf{0}
\end{bmatrix}, \qquad \mathbf{C} := (H\alpha)^T \otimes I_n,
\end{equation}and in view of (\ref{equ:rs_neuronal_model}),
we can write the linear DCM as a classical stochastic state-space model
\begin{equation}\label{equ:linear_stochastic_dcm}
\left\{\begin{array}{rcl}
\mathbf{x}(k+1) & = & \mathbf{A} \mathbf{x}(k) + \mathbf{w}(k) \\
y(k) &=& \mathbf{C} \mathbf{x}(k) + e(k).
\end{array}\right.
\end{equation}
Here, $\mathbf{w}(k):= [w^T(k) \ \mathbf{0}]^T\in\mathbb{R}^{ns}$ and  we assume
\begin{align}
\mathbf{w}(k) &\sim \mathcal{N}(\mathbf{0}, \mathbf{Q}), \qquad \mathbf{Q} :=\mathrm{blkdiag}(Q,\mathbf{0})\label{equ:process_noise}\\
e(k) &\sim\mathcal{N}(\mathbf{0},R), \qquad R:=\lambda^2 I_n\label{equ:measurement_noise}
\end{align}
where $\mathrm{blkdiag(\cdot)}$ is the block-diagonal operator and the covariance matrix $Q$ has been defined in \eqref{equ:Q}.
Notice that model \eqref{equ:linear_stochastic_dcm} only depends on the neuronal states $x\in\mathbb{R}^{n}$ and on the parameters $\eta=\{A, \alpha, \sigma,\lambda\}$. Next sections will describe how these are estimated, starting from the simpler (but unrealistic) situation in which the neuronal states $x$ are measurable, Sec.~\ref{sec:estimation_measured_x}, to the situation in which $x$ are not measurable, Sec.~\ref{sec:estimation_unmeasured_x}.

\section{Parameters Estimation: Measured Neuronal Activity}\label{sec:estimation_measured_x}
In this Section we unrealistically assume that observations of the neuronal activity (encoded in the states $x$) are available  and we detail how the parameters defining model \eqref{equ:linear_stochastic_dcm} are estimated.

\subsection{Estimation of matrix $A$}\label{sec:reweighted}

The matrix $A$ is expected to be sparse, encoding the fact that each $x_i(t)$ is influenced by the activities of  few other areas $x_j(t)$. To induce sparsity in the estimated $A$, we use a reweighted $\ell_1$ method, \cite{wipf2010iterative,wipf2008new}, originally proposed for linear regression models. However, because of the nonlinearity (in the parameters) of model \eqref{equ:rs_neuronal_model}, such approach cannot be straightforwardly applied but needs to be properly adapted.\\
Introducing the matrices
\begin{align}\label{equ:data_mat_A}
X_+ & :=\left[\begin{array}{c} x^T(2) \\ \vdots  \\ x^T(N) \end{array}\right], \; \;X:=\left[\begin{array}{c} x^T(1) \\ \vdots  \\ x^T(N-1) \end{array}\right],\nonumber \\W & :=\left[\begin{array}{c} w^T(1) \\ \vdots  \\ w^T(N-1) \end{array}\right].
\end{align}
we rewrite (\ref{equ:rs_neuronal_model}) in the non-linear regression form 
\begin{align} \label{sampled_model2}
X_+=X e^{A^T T_R} + W.
\end{align}
which can be linearised  exploiting  $e^{A^T T_R}\simeq I +A^T T_R$:
\begin{align}
\label{sampled_model_approx}\Delta X=XA^T T_R+W
\end{align}
where $\Delta X=X_+-X$. 
Using the vectorization operator, we can rewrite (\ref{sampled_model_approx}) as follows 
\begin{align}\label{sampled_model_approx2}
\mathrm{x}=\Phi a+\mathrm{w}
\end{align}
where $\mathrm{x}:=\mathrm{vec}(\Delta X)$, $\Phi:=(I\otimes  X) T_R$, $a:=\mathrm{vec}(A^T)$, $\mathrm{w}:=\mathrm{vec} (W)$. Adopting the Sparse Bayesian Learning perspective \cite{tipping2001sparse}, we model $a\sim \mathcal{N}(0,\Gamma)$ with $\Gamma:=\mathrm{diag}(\gamma_1 \ldots \gamma_{n^2})$ and independent of $\mathrm{w}$, so that  the covariance matrix of $\mathrm{x}$ is $\mathrm{Var}[\mathrm{x}]=\Phi \Gamma \Phi +Q\otimes I_N$. Direct application of the reweighted $\ell_1$ method  \cite{wipf2010iterative} to the approximated linear model (\ref{sampled_model_approx2}) is not possible here, since $T_R=2$s is far to be small, (\ref{sampled_model_approx2}) is only a crude approximation of the continuous time model. For this reason,  we propose Algorithm \ref{algo_Reweighted}
\begin{algorithm}
\caption{Proposed reweighted $\ell_1$ method}\label{algo_Reweighted}\begin{algorithmic}[1]
\Statex \textbf{Inputs:} $x(k), \ k=1,...,N$
\Statex \textbf{Initialization:} $A^{(0)}=-I_n$, $\gamma_i^{(0)}=0.25$ $i=1\ldots n^2$, $j=0$
\Repeat 
\State $\sigma^{(j)^2} = \frac{1}{Nn}\| X_+-X e^{A^{(j)^T} T_R}\|^2_{Q^{-1}}$
\State $Q= \sigma^{(j)^2} \int_0^{T_R} e^{A^{(j)}\tau} e^{A^{(j)^T}\tau}d\tau$
\State $A^{(j+1)}=\underset{A\in \mathcal A}{\mathrm{argmin}} \; \; \left\{\| X_+-X e^{A^T T_R}\|^2_{Q^{-1}}\right.$\label{alg_step:opt_probl_A}
\Statex \hspace{3cm}$ \left.+\mathrm{vec}(A^T)^T \left[\Gamma^{(j)}\right]^{-1} \mathrm{vec}(A^T)\right\}$
\State $\gamma_i^{(j+1)}=(a^{(j+1)}_i)^2+\gamma_i^{(j)}$
\Statex \hspace{2cm}$-(\gamma_i^{(j)})^2\phi_i^T( \Phi \Gamma^{(j)} \Phi^T +Q\otimes  I_N)^{-1}\phi_i $,
\Statex \hspace{.4cm} $i=1,...,n^2$
\State $k = k+1$
\Until{$\|A^{(j)}-A^{(j-1)}\|_F/\|A^{(j)}\|_F$ is sufficiently small}
\Statex \textbf{Outputs:} $A^{(j)}, \sigma^{(j)}$
\end{algorithmic}
\end{algorithm}
where $a^{(j)}=\mathrm{vec}(A^{(j)})$; $a_i^{(j)}$ is the entry in position $i$ in vector $a^{(j)}$; $\Gamma^{(j)}=\mathrm{diag}(\gamma_1^{(j)} \ldots \gamma_{n^2}^{(j)})$; $\phi_i$ is the $i$-th column of $\Phi$. Here the update for $\gamma_i$'s is the same as  \cite[formula (29)]{wipf2010iterative} while the main difference rests in 
update for  $A^{(j)}$, which  is set as the posterior mean of (\ref{sampled_model2}), i.e. using the exact model, rather than the posterior mean of the approximated model (\ref{sampled_model_approx2}); more precisely: 
\begin{equation}\label{equ:linear_solution_reweighted}
a^{(j+1)} = \left(\Phi^T Q^{-1}\Phi + [\Gamma^{(j)}]^{-1}\right)^{-1} \Phi^T Q^{-1}\mathrm{x}.
\end{equation}
Notice that, the search of the optimal $A^{(j)}$ is restricted to the set $\mathcal A$. In this paper, $\mathcal A$ is chosen as the set of stable matrices, i.e. with eigenvalues having negative real part. 

\begin{rem}
It is worth to stress that the proposed reweighted scheme automatically selects a sparsity pattern for the matrix $A$. This represents a significant computational savings w.r.t. the Variational Bayes methods mentioned in Sec.~\ref{sec:existing_methods}, which instead
have to be provided with several candidate sparsity patterns, among which the best is selected  resorting to model selection techniques, e.g. using the free-energy \cite{friston2011post}; this gives rise to a combinatorial explosion of candidates making practical inference infeasible. 
\end{rem}

\subsection{Estimation of the coefficients $\alpha$}
In Sec.~\ref{sec:linear_dcm}, a linear model for the hemodynamic response has been derived. As previously remarked, this depends on some parameters $\alpha$ that have to be estimated from data. Here we assume to have measured $N$ data-pairs $(x(k),y(k))$ related by the output equation of model \eqref{equ:linear_stochastic_dcm}, that is
\begin{equation}
y(k) =\mathbf{C}\mathbf{x}(k) +e(k), \qquad k=1,...,N.
\end{equation}
The previous equation can be rewritten as
\begin{equation}
y(k) =\varphi(k)H\alpha +e(k), \qquad k=1,...,N
\end{equation}
by defining $\mathrm{vec}(\varphi(k)):=\mathbf{x}(k)$ and observing  that $\mathbf{C}\mathbf{x}(k) = \mathrm{vec}(\varphi(k)H\alpha) = \varphi(k)H\alpha$. Consequently, introducing the matrices
\begin{equation}\label{equ:data_mat_alpha}
Y := \begin{bmatrix}
y(1)\\ \vdots \\ y(N) 
\end{bmatrix},  \quad \mathbf{H}:= \begin{bmatrix}
\varphi(1) \\ \vdots \\ \varphi(N)
\end{bmatrix}H, \quad E := \begin{bmatrix}
e(1)\\ \vdots \\ e(N)
\end{bmatrix}
\end{equation}
we have the following linear regression form
\begin{equation}
Y = \mathbf{H}\alpha + E.
\end{equation}
We consider  the Gaussian prior $\alpha\sim\mathcal{N}(\mu_\alpha, \Sigma_\alpha)$, where
\begin{equation}
\mu_\alpha := [1 \ 0 \ \cdots \ 0]^T, \qquad \Sigma_\alpha :=\mathrm{diag}(s_{1},...,s_{p})
\end{equation}
and $s_{i}$, $i=1,...,p$, are the first $p$ singular values of the empirical covariance matrix $\Sigma_h$ in equation \eqref{equ:sigma_h}. Recalling that $e(k)\sim\mathcal{N}(0,\lambda^2 I_n)$ and exploiting standard rules on conditioned Gaussian variables, we easily derive the MAP estimate of $\alpha$ as
\begin{equation}\label{equ:alpha_hat}
\hat{\alpha} = \mu_\alpha + (\mathbf{H}^T\mathbf{H} + \hat{\lambda}^2 \Sigma_\alpha^{-1})^{-1} \mathbf{H}^T (Y-\mathbf{H}\mu_\alpha)
\end{equation}
where
\begin{equation}\label{equ:lambda_hat}
\hat{\lambda}^2 = \frac{(Y-\mathbf{H}\hat{\alpha}_{LS})^T (Y-\mathbf{H}\hat{\alpha}_{LS})}{Nn-p-1}
\end{equation}
and $\hat{\alpha}_{LS} = (\mathbf{H}^T\mathbf{H})^{-1}\mathbf{H}^TY$.

\vspace{3mm}

\section{Parameters Estimation: Unmeasured Neuronal Activity}\label{sec:estimation_unmeasured_x}
This section deals with the realistic situation in which only the resting-state fMRI data are observed and the neuronal states $x$ are hidden. Sec.~\ref{sec:state_estimation} will describe how we estimate $x$ starting from the measured BOLD signal $y$, having fixed the model parameters $\eta=\{A, \alpha,\sigma,\lambda\}$; the subsequent Sec.~\ref{sec:params_estimation_hidden_x} proposes two algorithms for the combined estimation of $x$ and $\eta$.


\subsection{Neuronal Activity Estimation}\label{sec:state_estimation}
We assume that we are provided with $N$ observations of the BOLD signal $y$, modelled according to the linear dynamic model \eqref{equ:linear_stochastic_dcm}, specified by the parameters $\eta$ (here supposed to be known). We use the Rauch-Tung-Striebel smoother (RTSS) \cite{rauch1965maximum} to estimate the hidden neuronal states $x(k),\ k=1,...,N$, and in particular, the closed form smoothing solution to the linear model \eqref{equ:linear_stochastic_dcm}, that is
\begin{equation}\label{equ:smoothing_solution}
p(\mathbf{x}(k) | Y,\ \eta) = \mathcal{N}(\hat{\mathbf{x}}^s(k), \mathbf{P}^s(k)).
\end{equation}
The implementation of the RTSS is summarized in Algorithm \ref{alg:smoother}.

\begin{algorithm}
\caption{RTS Smoother}\label{alg:smoother}
\begin{algorithmic}[1]
\Statex \textbf{Inputs:} $y(k),\ k=1,..,N$;  $\mathbf{A}, \mathbf{C}$ in Eq. \eqref{equ:a_c}, $\mathbf{Q}$ and $R$ in Eqs. \eqref{equ:process_noise}, \eqref{equ:measurement_noise}.
\Statex \textbf{Forward Recursion}
\State Initialize: $\hat{\mathbf{x}}(0)=0$, $\mathbf{P}(0)=I_{ns}$
\For{$k=1,...,N$}
\State $\hat{\mathbf{x}}^-(k) = \mathbf{A}\hat{\mathbf{x}}(k-1)$
\State $\mathbf{P}^-(k) = \mathbf{A}\mathbf{P}(k-1)\mathbf{A}^T + \mathbf{Q}$
\State $S(k) = \mathbf{C} \mathbf{P}^-(k) \mathbf{C} + R$
\State $K(k) = \mathbf{P}^-(k) \mathbf{C}^T S^{-1}(k)$
\State $\hat{\mathbf{x}}(k) = \hat{\mathbf{x}}^-(k) + K(k) [y(k) -\mathbf{C} \hat{\mathbf{x}}^-(k)]$
\State $\mathbf{P}(k) = \mathbf{P}^-(k) - K(k) S(k) K^T(k)$
\EndFor
\Statex \textbf{Backward recursion}
\State Initialize: $\hat{\mathbf{x}}^s(N) = \hat{\mathbf{x}}(N), \ \mathbf{P}^s(N)=\mathbf{P}(N)$
\For{$k=N-1,...,0$}
\State $\hat{\mathbf{x}}^-(k+1) = \mathbf{A} \hat{\mathbf{x}}(k)$
\State $\mathbf{P}^-(k+1) = \mathbf{A} \mathbf{P}(k) \mathbf{A}^T + \mathbf{Q}$
\State $\mathbf{G}(k)=\mathbf{P}(k) \mathbf{A}^T [\mathbf{P}^-(k+1)]^{-1}$
\State $\hat{\mathbf{x}}^s(k)=\hat{\mathbf{x}}(k) + \mathbf{G}(k) [\hat{\mathbf{x}}^s(k+1) - \hat{\mathbf{x}}^-(k+1)]$
\State $\mathbf{P}^s(k) = \mathbf{P}(k) + \mathbf{G}(k) [\mathbf{P}^s(k+1) - \mathbf{P}^-(k+1)] \mathbf{G}^T(k)$
\EndFor
\Statex \textbf{Outputs:} $\hat{\mathbf{x}}^s(k),\ \mathbf{P}^s(k),\ \mathbf{G}(k), \ k=1,...,N$
\end{algorithmic}
\end{algorithm}

\subsection{Parameters Estimation}\label{sec:params_estimation_hidden_x}
We turn to the problem of estimating the parameters $\eta$ specifying the linear model \eqref{equ:linear_stochastic_dcm}. Adopting the Bayesian perspective, we would like to determine the value of $\eta$ which maximizes the marginal posterior
\begin{equation}\label{equ:posterior_eta}
p(\eta |Y) = \int p(\mathbf{X}, \eta | Y)\ \mathrm{d}\mathbf{X}
\end{equation}
where $\mathbf{X}:=[\mathbf{x}^T(0) \ \cdots \ \mathbf{x}^T(N)]^T$. However, the computation of such a high-dimensional integral is typically avoided by exploiting the decomposition $p(\eta|Y) \propto p(Y|\eta)p(\eta)$ and resorting to an iterative algorithm to compute the likelihood
\begin{equation}
p(Y|\eta) = \int p(\mathbf{X},Y|\eta) \ \mathrm{d}\mathbf{X}.
\end{equation}
Here we assume $p(\eta)\propto p(A)p(\alpha)p(\sigma)p(\lambda)$, where $p(\sigma)$ and $p(\lambda)$ are uninformative priors, while $p(\mathrm{vec}(A^T))\sim \mathcal{N}(\mathbf{0},\Gamma)$ and $p(\alpha)\sim\mathcal{N}(\mu_\alpha,\Sigma_\alpha)$.
\\In this work we use the Expectation-Maximization (EM) algorithm to maximize $\ln p(Y|\eta)$, which is equivalent to maximizing $p(Y| \eta)$. The EM method performs such optimization by iteratively maximizing the lower bound
 \begin{equation}\label{equ:em_lower_bound}
 \mathcal{F} (q(\mathbf{X}),\eta) = \int q(\mathbf{X}) \ln \frac{p(\mathbf{X},Y|\eta)}{q(\mathbf{X})}\ \mathrm{d}\mathbf{X} 
 \end{equation}
 of $\ln p(Y|\eta)$,
 w.r.t. an arbitrary distribution $q(\mathbf{X})$ and $\eta$. In the statistical learning literature $ \mathcal{F} (q(\mathbf{X}),\eta) $ is also known as \textit{free-energy}.
 At the $l$-th iteration of the algorithm, $\mathcal{F} (q(\mathbf{X}),\eta^{(l)})$ is maximized by $q^{(l+1)}(\mathbf{X})= p(\mathbf{X}|Y,\eta^{(l)})$. Plugging this into \eqref{equ:em_lower_bound}, one obtains
 \begin{align}
\mathcal{F}(q^{(l+1)}(\mathbf{X}),\eta)  = &\int p(\mathbf{X}|Y,\eta^{(l)})\ln p(\mathbf{X},Y|\eta)\ \mathrm{d}\mathbf{X} \nonumber\\
&-\int p(\mathbf{X}|Y,\eta^{(l)})\ln p(\mathbf{X}|Y\eta^{(l)})\ \mathrm{d}\mathbf{X} 
 \end{align}
which now needs to be maximized w.r.t. $\eta$. Noticing that the second term of $\mathcal{F}(q^{(l+1)}(\mathbf{X}),\eta)$ does not depend on $\eta$, this coincides with optimizing
\begin{equation}\label{equ:em_Q}
\mathcal{Q}(\eta,\eta^{(l)}) = \int p(\mathbf{X}|Y,\eta^{(l)})\ln p(\mathbf{X},Y|\eta) \ \mathrm{d}\mathbf{X}
\end{equation}
w.r.t. $\eta$. Using the Markovian property of system \eqref{equ:linear_stochastic_dcm}, $\mathcal{Q}(\eta,\eta^{(l)})$ can be rewritten as (\cite{sarkka2013bayesian}, Ch.12)
\begin{align}
\mathcal{Q}(\eta,\eta^{(l)})= &\sum_{k=1}^N\int p(\mathbf{x}(k), \mathbf{x}(k-1)|Y,\eta^{(l)}) \times \label{equ:em_Q2}\\
& \qquad\quad \times\ln p(\mathbf{x}(k)|\mathbf{x}(k-1),\eta)\mathrm{d}\mathbf{x}(k)\mathrm{d}\mathbf{x}(k-1)\nonumber\\
&+ \sum_{k=1}^N\int p(\mathbf{x}(k)|Y,\eta^{(l)})\ln p(y(k)|\mathbf{x}(k),\eta)\mathrm{d}\mathbf{x}(k).\nonumber
\end{align} 
Plugging the smoothing distributions $p(\mathbf{x}(k)|Y,\eta^{(l)})= \mathcal{N}(\hat{\mathbf{x}}^s(k), \mathbf{P}^s(k))$ (see Eq. \eqref{equ:smoothing_solution}) and the pairwise smoothing distributions (\cite{sarkka2013bayesian}, Ch.12)
\begin{align}
p&(\mathbf{x}(k),\mathbf{x}(k-1)|Y,\eta^{(l)})=\\
& \mathcal{N} \left(
\begin{bmatrix}
\hat{\mathbf{x}}^s(k)\\\hat{\mathbf{x}}^s(k-1)
\end{bmatrix}, \begin{bmatrix}
\mathbf{P}^s(k) & \mathbf{P}^s(k) \mathbf{G}^T(k-1)\\ \mathbf{G}(k-1)\mathbf{P}^s(k) & \mathbf{P}^s(k-1)
\end{bmatrix}
 \right)\nonumber
\end{align}
into \eqref{equ:em_Q2} we get
\begin{align}
\mathcal{Q}(\eta,\eta^{(l)})= &-\frac{N}{2} \ln|2\pi\mathbf{Q}| -\frac{N}{2} \ln|2\pi R| \label{equ:em_Q3}\\
&-\frac{N}{2}\mbox{tr}\left[ \mathbf{Q}^{-1}\left(\Lambda - \Psi \mathbf{A}^T -\mathbf{A}\Psi^T + \mathbf{A}\Upsilon\mathbf{A}^T \right) \right]\nonumber\\
&-\frac{N}{2}\mbox{tr}\left[ R^{-1}\left(\Delta - \Xi \mathbf{C}^T -\mathbf{C}\Xi^T + \mathbf{C}\Lambda \mathbf{C}^T \right) \right]\nonumber
\end{align}
where
\begin{align}
\Lambda &= \frac{1}{N}\sum_{k=1}^N \mathbf{P}^s(k) + \hat{\mathbf{x}}^s(k)\left[\hat{\mathbf{x}}^s(k)\right]^T\nonumber\\
\Upsilon &= \frac{1}{N}\sum_{k=1}^N \mathbf{P}^s(k-1) + \hat{\mathbf{x}}^s(k-1)\left[\hat{\mathbf{x}}^s(k-1)\right]^T\nonumber\\
\Xi &= \frac{1}{N}\sum_{k=1}^N  y(k)\left[\hat{\mathbf{x}}^s(k)\right]^T\nonumber\\
\Psi &= \frac{1}{N}\sum_{k=1}^N \mathbf{P}^s(k) \mathbf{G}(k-1)+ \hat{\mathbf{x}}^s(k)\left[\hat{\mathbf{x}}^s(k-1)\right]^T\nonumber\\
\Delta &= \frac{1}{N}\sum_{k=1}^N y(k)y^T(k).\nonumber
\end{align}
The complete routine we adopt to estimate the parameters $\eta$ through the EM method is summarized in Algorithm \ref{alg:em}. The procedure is initialized with $\eta^{(0)} = \{A^{(0)}, \sigma^{(0)}, \alpha^{(0)}, \lambda^{(0)}\}$:
\begin{align}\label{equ:eta_init}
A^{(0)} &=-I_n, \qquad \sigma^{(0)}=10^{-2}\\
\alpha^{(0)}&=\mu_\alpha, \qquad \  \lambda^{(0)} = \sqrt{\frac{
\mathrm{tr}[(Y-\bar{Y})(Y-\bar{Y})^\top]}{10Nn}} \nonumber
\end{align}
where $\bar{Y}=\frac{1}{Nn}\sum_{i=1}^{Nn} Y_i$.

\begin{algorithm}
\caption{Estimation of parameters $\eta$ through EM}\label{alg:em}
\begin{algorithmic}[1]
\Statex \textbf{Inputs:} $y(k),\ k=1,..,N$, $H$, $\mu_\alpha$, $\Sigma_\alpha$
\Statex \textbf{Initialization:}  Set $\eta^{(0)}$ as in Eq. \eqref{equ:eta_init}, set $l=0$
\Repeat
\State Apply Algorithm \ref{alg:smoother} to get $\hat{\mathbf{x}}^s(k)$, $\mathbf{P}^s(k)$, $\mathbf{G}(k)$,
\Statex \hspace{.6cm}$k=1,...,N$
\State Compute $\mathcal{Q}(\eta, \eta^{(l)})$ using Eq. \eqref{equ:em_Q3}
\State $\eta^{(l+1)} = \arg\max_{\eta\in \Omega} \left\{\mathcal{Q}(\eta, \eta^{(l)})\right.$\label{alg_step:Q_max}

\vspace{1.5mm}

\Statex \hspace{3cm}$- \frac{1}{2}\ln|2\pi\Gamma^{(l)}| -\frac{1}{2}\ln|2\pi\Sigma_\alpha|$

\vspace{1.5mm}

\Statex\hspace{3cm}$- \frac{1}{2}\mathrm{vec}(A^T)^T [\Gamma^{(l)}]^{-1}\mathrm{vec}(A^T) $

\vspace{1.5mm}

\Statex \hspace{3cm}$\left. -\frac{1}{2}(\alpha - \mu_\alpha)\Sigma_\alpha^{-1}(\alpha - \mu_\alpha)\right\}$

\vspace{1.5mm}

\State $Q= \sigma^{(l+1)^2} \int_0^{T_R} e^{A^{(l+1)}\tau} e^{A^{(l+1)^T}\tau}d\tau$

\vspace{1.5mm}

\State $\gamma_i^{(l+1)}=(a^{(l+1)}_i)^2+\gamma_i^{(l)}$
\Statex \hspace{2cm}$-(\gamma_i^{(l)})^2\phi_i^T( \Phi \Gamma^{(l)} \Phi^T +Q\otimes  I_N)^{-1}\phi_i $,
\Statex \hspace{.4cm} $i=1,...,n^2$
\State $l=l+1$
\Until{$\|A^{(l)}-A^{(l-1)}\|_F/\|A^{(l)}\|_F$ is sufficiently small}
\Statex \textbf{Outputs:} $\eta^{(l)}$
\end{algorithmic}
\end{algorithm}

\noindent Notice that, at Step \ref{alg_step:Q_max} of Algorithm \ref{alg:em}, the search of $\eta$ is restricted to the set $\Omega=\{\mathcal{A},\ \mathbb{R}_+,\ \mathbb{R}^{p+1},\ \mathbb{R}_+ \}$. Since this optimization may be computationally intensive, we have also tested a simpler method, where  the uncertainty of the estimates returned by the smoother (encoded in the matrix $\mathbf{P}^s(k)$) is ignored. Briefly, this alternative routine  iteratively estimates $A$ and $\sigma$ through Algorithm \ref{algo_Reweighted}, $\alpha$ and $\lambda$ through Eqs. \eqref{equ:alpha_hat}, \eqref{equ:lambda_hat}, replacing the hidden neuronal states with the smoother estimates $\hat{\mathbf{x}}^s(k)$. The results are inferior to those of Algorithm \ref{alg:em} and thus not reported. 

\section{Experimental Results}\label{sec:experiments}
The methods illustrated in Secs.~\ref{sec:estimation_measured_x} and \ref{sec:estimation_unmeasured_x} are experimentally evaluated using simulated rs-fMRI data. In particular: (a) we analyze the effectiveness of the  reweighted $\ell_1$ technique detailed in Sec.~\ref{sec:reweighted} in detecting the correct brain network structure and in estimating the strength of the neuronal connections; (b) we evaluate the performance of  Algorithm \ref{alg:em} which only uses  the  resting state BOLD signals, when the true network structure has to be inferred from the fMRI data. In the latter setup, the proposed methods are compared with the Variational Bayes routines implemented in the SPM12 toolbox (\texttt{http://www.fil.ion.ucl.ac.uk/spm/}).

\subsection{Experimental Setup}\label{sec:experimental_setup}
The rs fMRI data used in these experiments are generated through the SPM12 routines \verb!spm_int_J!, \verb!spm_fx_fmri! and \verb!spm_gx_fmri!. Specifically, we consider a brain network with $n=7$ neuronal nodes, whose activity is modulated according to Eq. \eqref{equ:stochastic_bilinear_neuronal_model}, where $w(t)$ is a Gaussian white noise process with variance $\sigma^2=0.01$ and 
\begin{equation}\label{equ:true_A}
\footnotesize A = \begin{bmatrix}
-0.5 & 0 & 0 & 0 & -0.2 & 0 & 0\\
0 & -0.5 & 0 & -0.45 & -0.3 & 0 & 0\\
0 & 0 & -0.5 & 0.8 & 0 & 0 & 0\\
0 & 0.6 & 0 & -0.5 & -0.1 & 0.6 &0\\
0.3 & 0 & -0.55 & 0 & -0.5 & 0.2 & 0\\
0 & 0 & 0 & 0 & 0.3 & -0.5 & 0.45\\
0.15 & 0 & 0.2 & 0 & 0 & 0 & -0.5
\end{bmatrix}
\end{equation} 
$N=600$ data are generated with sampling time $T_R=2$s. No observation noise 
is added to the BOLD signal generated by the routine \verb!spm_gx_fmri!.  It is worth stressing that the self-connections are a-priori fixed to $-0.5$ to prevent instability issues.\\

\subsection{Performance Evaluation}
Let $\hat{A}$ be an estimate of the effective connectivity matrix: we define the 
Root-Mean-Squared-Error (RMSE) 
\begin{equation}\label{equ:rmse}
RMSE(\hat{A}) = \frac{\|\underline{A} - \underline{\hat{A}}\|_F}{\sqrt{n(n-1)}}
\end{equation}
where the notation $\underline{A}$ denotes the matrix $A$ with its diagonal set to 0. Indeed, the diagonal elements of $A$ (i.e. the self-connections)  do not give any information on the effective connectivity.\\
A second performance index  measures the number of errors committed by a certain method in the recovery of the true sparsity pattern of matrix \eqref{equ:true_A}. 
 Namely,
\begin{equation}\label{equ:err}
ERR(\hat{A}) = \|SP(A)-SP(\hat{A})\|_F^2 
\end{equation}
where $SP(\cdot): \mathbb{R}^{n\times m} \rightarrow \mathbb{R}^{n\times m}$ and
\begin{equation}
[SP(A)]_{ij}:=\left\{\begin{array}{ll}
1 & \mbox{if } [A]_{ij}\neq 0\\
0 & \mbox{if } [A]_{ij} = 0
\end{array}\right..
\end{equation}

\subsection{Estimation of effective connectivity using measured neuronal activity}
In this section we assess the validity of the reweighted $\ell_1$ method illustrated in Sec.~\ref{sec:reweighted}. To this purpose we pass as inputs to Algorithm \ref{algo_Reweighted} the neuronal time-series $x(k)$, $k=1,...,600$, generated by the SPM12 routine \verb!spm_int_J! when provided with \verb!spm_fx_fmri! and the matrix \eqref{equ:true_A}. In addition, we compare Algorithm \ref{algo_Reweighted} with the linear  reweighted $\ell_1$ based on approximation \eqref{sampled_model_approx2}.\\
Table \ref{tab:reweighted_results} compares the results achieved by Algorithm \ref{algo_Reweighted} and by its linear version on the data generated in 50 Monte-Carlo runs, using different realizations of the noise $w(t)$: the medians (and the standard deviations) of indexes \eqref{equ:rmse} and \eqref{equ:err} are reported. Before evaluating these indexes, a tresholding is applied to the estimated matrices $\hat{A}$: namely the entries with absolute value less than 0.1 are set to zero.\\
The results reported in Table \ref{tab:reweighted_results} show that Algorithm \ref{algo_Reweighted} benefits from using the true generative model at Step \ref{alg_step:opt_probl_A}, in place of the linear approximation \eqref{sampled_model_approx2}.

\begin{table}
\centering
\begin{tabular}{l|cc}
\toprule
&  $ERR(\hat{A})$ & $RMSE(\hat{A})$ \\
\midrule
Algorithm \ref{algo_Reweighted} & 3 (2.9) & 0.05 (0.11)  \\
Linear reweighted $\ell_1$ algorithm &  10 (1.4) & 0.18 (0.006)\\
\bottomrule
\end{tabular}
\caption{Performance of Algorithm \ref{algo_Reweighted} and classical reweighted $\ell_1$ method over 50 Monte-Carlo runs (using measured neuronal activity $x(k),\ k=1,..,600$. Median and std (between brackets) of metrics \eqref{equ:rmse} and \eqref{equ:err} are reported.}\label{tab:reweighted_results}
\end{table}

{

\subsection{Estimation of effective connectivity using rs fMRI data}
We now evaluate the performance of the proposed Algorithm \ref{alg:em}  in the recovery of matrix \eqref{equ:true_A} when the BOLD data $y(k)$, $k=1,...,600$, generated under the setup of Sec.~\ref{sec:experimental_setup}, are available. In particular, we conduct a Monte-Carlo study over 20 different realizations of the noise $w(k)$. The performance of Algorithm \ref{alg:em} is compared with the one returned by the SPM12 routines \verb!spm_dcm_fmri_csd! and \verb!spm_dcm_estimate!, both provided with null exogenous inputs. These functions implement two of the methods mentioned in Sec.~\ref{sec:existing_methods}: specifically, \verb!spm_dcm_fmri_csd! runs the spectral DCM (\textbf{sDCM}) \cite{friston2014dcm}, while \verb!spm_dcm_estimate! is invoked activating the option \verb!DCM.options.stochastic! in order to implement the \textbf{DEM} procedure \cite{friston2008variational}.\\
To guarantee a fair comparison between our methods and the SPM12 routines, the latter should perform model selection among all the possible sparsity patterns of matrix $A$. However, this would require the estimation of a combinatorial number of matrices, making this type of approach computationally infeasible. Attempting to still conduct a fair comparison, we pass 21 candidate sparsity patterns $\mathsf{m}_i$, $i=1,...,21$, to the SPM12 routines, including the true one (denoted by \textbf{(a)}). These are illustrated in Figure \ref{fig:sparsity_patterns}. For each of them, Table \ref{tab:friston_results_12} reports the value of index \eqref{equ:err}, as well as the median and the standard deviation of index \eqref{equ:rmse} corresponding to the matrices estimated by the routines \verb!spm_dcm_fmri_csd! and \verb!spm_dcm_estimate!. In addition, the columns denoted with $\#Chosen$ report the number of times each model is chosen according to the free-energy metric. The most frequently selected structures are very different from the true one (as indicated by the value of the metrics $ERR(\hat{A})$). According to our experimental experience, the SPM routines tend to select structures with a low degree of sparsity.\\
Table \ref{tab:our_alg_results} contains the median and the standard deviation of indexes \eqref{equ:rmse} and \eqref{equ:err} computed for the matrices $\hat{A}$ estimated by the proposed Algorithm \ref{alg:em}. The performance of our approach is superior to the one achieved by the SPM12 routines, both in terms of $RMSE(\hat{A})$ and $ERR(\hat{A})$.   In addition, the sparsity patterns detected by Algorithm \ref{alg:em} differ at most by 7 entries from the true one. It is worth noting that $RMSE(\hat{A})$ returned by Algorithm \ref{alg:em} is comparable to that achieved by spectral DCM when the latter is provided with the true sparsity pattern.\\
On the conducted study, Alg. \ref{alg:em} takes an average computational time equal to 558s, while \verb!spm_dcm_fmri_csd! and \verb!spm_dcm_estimate! respectively  take 
on average 721s and 4044s to process all the 21 patterns.

\begin{figure}
\centering
\includegraphics[width=\columnwidth]{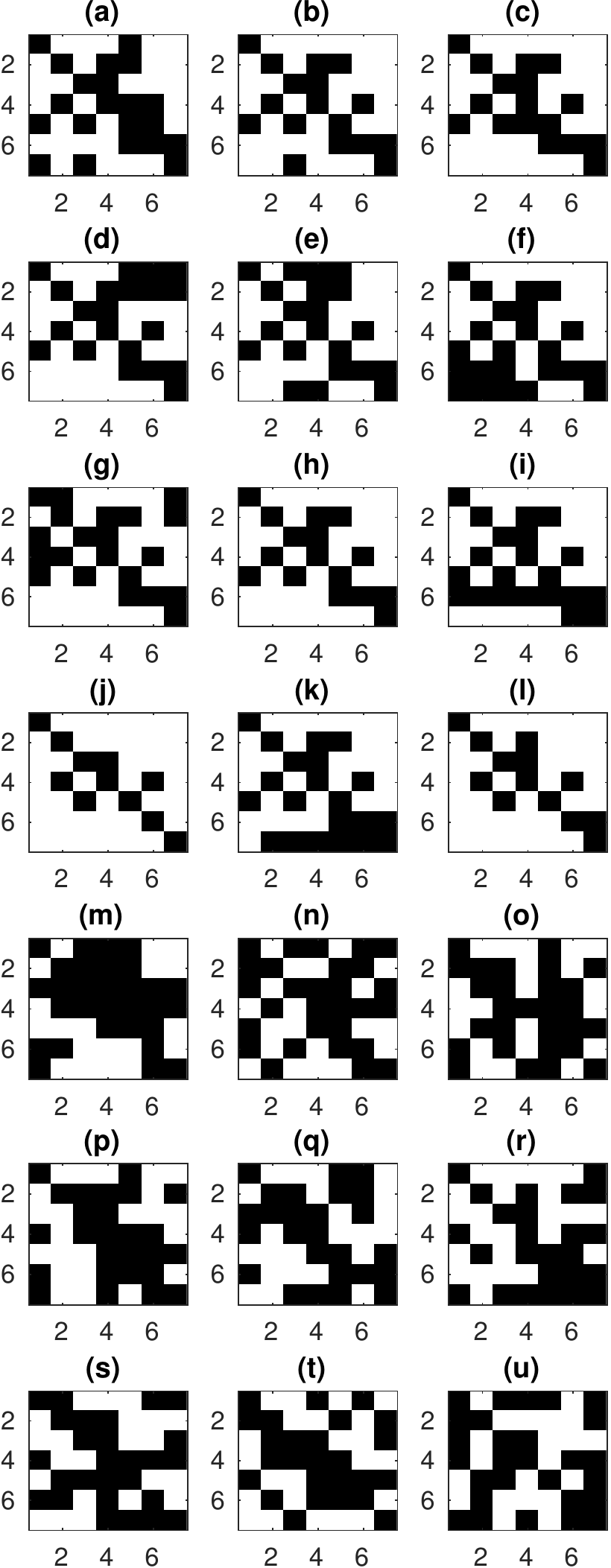}
\cprotect\caption{Candidate sparsity patterns $\mathsf{m}_i$ which are provided to the SPM12 routines. Black cells denote non-zero elements. \textbf{(a)} is the true one.}\label{fig:sparsity_patterns}
\end{figure}

\begin{table}
\centering
\begin{tabular}{l|c|cc|cc}
\toprule
& &\multicolumn{2}{c|}{\textbf{spectral DCM}} & \multicolumn{2}{c}{\textbf{DEM}} \\
&$ERR(\hat{A})$ & $\#Chosen$ & $RMSE(\hat{A})$& $\#Chosen$  & $RMSE(\hat{A})$ \\ 
\midrule
\textbf{(a)} & 0& 2 & 0.12 (0.05) & 0 & 0.21 (0.004)\\
\textbf{(b)} & 4& 0 & 0.25 (0.02)  & 0 & 0.21 (0.005)\\
\textbf{(c)} &6 & 0 & 0.23 (0.03)  & 0 & 0.21 (0.003)\\
\textbf{(d)} &8 & 0 & 0.18 (0.02) &   0 & 0.21 (0.004)\\
\textbf{(e)} & 6 & 1 & 0.14 (0.05) &  0 & 0.21 (0.005)\\
\textbf{(f)} & 8& 0& 0.22 (0.04)  & 0 & 0.21 (0.005)\\
\textbf{(g)} & 10& 0& 0.28 (0.09) & 0 & 0.22 (0.007)\\
\textbf{(h)} & 5& 0& 0.22 (0.02) & 0 & 0.20 (0.003)\\
\textbf{(i)} & 10& 0&  0.24 (0.07) & 0 & 0.21 (0.004)\\
\textbf{(j)} & 10 & 0& 0.17 (0.004)  & 0 & 0.22 (0.002)\\
\textbf{(k)} & 8&0 & 0.27 (0.03) & 0 & 0.21 (0.005)\\
\textbf{(l)} & 8& 0 & 0.18 (0.01) & 0 & 0.22 (0.002)\\
\textbf{(m)} & 19& 1 & 0.33 (0.04) & 3 & 0.23 (0.003)\\
\textbf{(n)} & \textbf{22 }& \textbf{14} & \textbf{0.31 (0.06)} & 2 & 0.23 (0.005)\\
\textbf{(o)} & \textbf{19} & 0 & 0.36 (0.04) & \textbf{14} & \textbf{0.23 (0.003)}\\
\textbf{(p)} & 15& 0 & 0.33 (0.04) & 0 & 0.22 (0.002)\\
\textbf{(q)} & 21& 2 & 0.32 (0.04) & 0 & 0.23 (0.003)\\
\textbf{(r)} & 16& 0 & 0.28 (0.02) & 0 & 0.22 (0.002)\\
\textbf{(s)} & 24& 0 & 0.31 (0.03) & 0 & 0.23 (0.004)\\
\textbf{(t)} & 17& 0 & 0.28 (0.03) & 1 & 0.23 (0.002)\\
\textbf{(u)} & 23& 0 & 0.27 (0.02) & 0 & 0.22 (0.002)\\
\bottomrule
\end{tabular}
\cprotect\caption{Performance of the SPM12 routines \verb|spm_dcm_fmri_csd| (\textbf{sDCM}) and  \verb|spm_dcm_estimate| (\textbf{DEM}) over 20 Monte-Carlo runs,
when they are given the 21 sparsity patterns of Fig. \ref{fig:sparsity_patterns}.
Median and std (within brackets) of index \eqref{equ:rmse} are reported, together with the no. of times ($\#Chosen$) each sparsity pattern is chosen according to the free-energy metric. Bold font highlights the most frequently selected models.}\label{tab:friston_results_12}
\end{table}

\begin{table}
\centering
\begin{tabular}{l|ccc}
\toprule
& $ERR(\hat{A})$ &  $RMSE(\hat{A})$  \\
\midrule
Algorithm \ref{alg:em} & 4 (0.74) & 0.13 (0.008) \\
\bottomrule
\end{tabular}\caption{Performance of Algorithm \ref{alg:em} over 20 Monte-Carlo runs. Median and std (within brackets) of the metrics \eqref{equ:rmse} and \eqref{equ:err} are reported.}\label{tab:our_alg_results}
\end{table}
}

\section{Conclusions and Future Work}\label{sec:conlusion}
We have proposed a new approach for the estimation of the effective connectivity in brain networks using resting state fMRI data. Existing methods tackle the problem by formulating a parametric generative model (DCM) of the observed data and by inverting it using Bayesian approaches. However, due to the nonlinear and stochastic nature of the standard DCMs, the current techniques are highly computationally intensive. In addition, when no a-priori knowledge is available on the   structure of the underlying network connecting the monitored brain regions, these approaches require to compare several candidate interaction patterns in order to choose the most plausible one. \\
The technique we have developed is based on a linear stochastic generative model (DCM), which we have built by linearly approximating the Balloon model describing the hemodynamic response. Exploiting the linearity of the derived DCM we have applied classical smoothing techniques to estimate both the hidden neuronal activity and the model parameters starting from the observed BOLD signal. A key feature of our algorithms is the ability to automatically detect the structure of the underlying brain network, thanks to the use of a sparsity inducing   method, thus avoiding a combinatorial search over candidate structures.\\
Future work will include an extended experimental study of the proposed method and, in particular, its application to real fMRI data.


\begin{thebibliography}{10}
\providecommand{\url}[1]{#1}
\csname url@samestyle\endcsname
\providecommand{\newblock}{\relax}
\providecommand{\bibinfo}[2]{#2}
\providecommand{\BIBentrySTDinterwordspacing}{\spaceskip=0pt\relax}
\providecommand{\BIBentryALTinterwordstretchfactor}{4}
\providecommand{\BIBentryALTinterwordspacing}{\spaceskip=\fontdimen2\font plus
\BIBentryALTinterwordstretchfactor\fontdimen3\font minus
  \fontdimen4\font\relax}
\providecommand{\BIBforeignlanguage}[2]{{%
\expandafter\ifx\csname l@#1\endcsname\relax
\typeout{** WARNING: IEEEtran.bst: No hyphenation pattern has been}%
\typeout{** loaded for the language `#1'. Using the pattern for}%
\typeout{** the default language instead.}%
\else
\language=\csname l@#1\endcsname
\fi
#2}}
\providecommand{\BIBdecl}{\relax}
\BIBdecl

\bibitem{friston2011functional}
K.~J. Friston, ``Functional and effective connectivity: a review,'' \emph{Brain
  connectivity}, vol.~1, no.~1, pp. 13--36, 2011.

\bibitem{friston2003dynamic}
K.~J. Friston, L.~Harrison, and W.~Penny, ``Dynamic causal modelling,''
  \emph{Neuroimage}, vol.~19, no.~4, pp. 1273--1302, 2003.

\bibitem{marreiros2008dynamic}
A.~C. Marreiros, S.~J. Kiebel, and K.~J. Friston, ``Dynamic causal modelling
  for fmri: a two-state model,'' \emph{Neuroimage}, vol.~39, no.~1, pp.
  269--278, 2008.

\bibitem{stephan2008nonlinear}
K.~E. Stephan, L.~Kasper, L.~M. Harrison, J.~Daunizeau, H.~E. den Ouden,
  M.~Breakspear, and K.~J. Friston, ``Nonlinear dynamic causal models for
  fmri,'' \emph{Neuroimage}, vol.~42, no.~2, pp. 649--662, 2008.

\bibitem{li2011generalised}
B.~Li, J.~Daunizeau, K.~E. Stephan, W.~Penny, D.~Hu, and K.~Friston,
  ``Generalised filtering and stochastic dcm for fmri,'' \emph{Neuroimage},
  vol.~58, no.~2, pp. 442--457, 2011.

\bibitem{friston2014dcm}
K.~J. Friston, J.~Kahan, B.~Biswal, and A.~Razi, ``A dcm for resting state
  fmri,'' \emph{Neuroimage}, vol.~94, pp. 396--407, 2014.

\bibitem{friston2007variational}
K.~Friston, J.~Mattout, N.~Trujillo-Barreto, J.~Ashburner, and W.~Penny,
  ``Variational free energy and the laplace approximation,'' \emph{Neuroimage},
  vol.~34, no.~1, pp. 220--234, 2007.

\bibitem{friston2008variational}
K.~J. Friston, N.~Trujillo-Barreto, and J.~Daunizeau, ``Dem: a variational
  treatment of dynamic systems,'' \emph{Neuroimage}, vol.~41, no.~3, pp.
  849--885, 2008.

\bibitem{daunizeau2009variational}
J.~Daunizeau, K.~J. Friston, and S.~J. Kiebel, ``Variational bayesian
  identification and prediction of stochastic nonlinear dynamic causal
  models,'' \emph{Physica D: Nonlinear Phenomena}, vol. 238, no.~21, pp.
  2089--2118, 2009.

\bibitem{friston2011post}
K.~Friston and W.~Penny, ``Post hoc bayesian model selection,''
  \emph{Neuroimage}, vol.~56, no.~4, pp. 2089--2099, 2011.

\bibitem{buxton1998dynamics}
R.~B. Buxton, E.~C. Wong, and L.~R. Frank, ``Dynamics of blood flow and
  oxygenation changes during brain activation: the balloon model,''
  \emph{Magnetic resonance in medicine}, vol.~39, no.~6, pp. 855--864, 1998.

\bibitem{friston2000nonlinear}
K.~J. Friston, A.~Mechelli, R.~Turner, and C.~J. Price, ``Nonlinear responses
  in fmri: the balloon model, volterra kernels, and other hemodynamics,''
  \emph{NeuroImage}, vol.~12, no.~4, pp. 466--477, 2000.

\bibitem{stephan2007comparing}
K.~E. Stephan, N.~Weiskopf, P.~M. Drysdale, P.~A. Robinson, and K.~J. Friston,
  ``Comparing hemodynamic models with dcm,'' \emph{Neuroimage}, vol.~38, no.~3,
  pp. 387--401, 2007.

\bibitem{friston2010generalised}
K.~Friston, K.~Stephan, B.~Li, and J.~Daunizeau, ``Generalised filtering,''
  \emph{Mathematical Problems in Engineering}, vol. 2010, 2010.

\bibitem{razi2015construct}
A.~Razi, J.~Kahan, G.~Rees, and K.~J. Friston, ``Construct validation of a dcm
  for resting state fmri,'' \emph{Neuroimage}, vol. 106, pp. 1--14, 2015.

\bibitem{frassle2017regression}
S.~Fr{\"a}ssle, E.~I. Lomakina, A.~Razi, K.~J. Friston, J.~M. Buhmann, and
  K.~E. Stephan, ``Regression dcm for fmri,'' \emph{NeuroImage}, 2017.

\bibitem{wipf2010iterative}
D.~Wipf and S.~Nagarajan, ``Iterative reweighted $\ell\_1$ and $\ell_2 $
  methods for finding sparse solutions,'' \emph{IEEE Journal of Selected Topics
  in Signal Processing}, vol.~4, no.~2, pp. 317--329, 2010.

\bibitem{wipf2008new}
D.~P. Wipf and S.~S. Nagarajan, ``A new view of automatic relevance
  determination,'' in \emph{Advances in neural information processing systems},
  2008, pp. 1625--1632.

\bibitem{tipping2001sparse}
M.~E. Tipping, ``Sparse bayesian learning and the relevance vector machine,''
  \emph{Journal of machine learning research}, vol.~1, no. Jun, pp. 211--244,
  2001.

\bibitem{rauch1965maximum}
H.~E. Rauch, C.~Striebel, and F.~Tung, ``Maximum likelihood estimates of linear
  dynamic systems,'' \emph{AIAA journal}, vol.~3, no.~8, pp. 1445--1450, 1965.

\bibitem{sarkka2013bayesian}
S.~S{\"a}rkk{\"a}, \emph{Bayesian filtering and smoothing}.\hskip 1em plus
  0.5em minus 0.4em\relax Cambridge University Press, 2013, vol.~3.

\end{thebibliography}
\end{document}